\documentclass[11pt,fleqn]{article}
\usepackage{graphicx}
\usepackage{amsmath}
\usepackage{amsthm}
\usepackage{amssymb}
\usepackage{amsfonts}
\usepackage{epsfig}
\usepackage{color}

\def\arg {\mathop{\rm arg}\nolimits}

\def\cot {{\rm cot}}

\newtheorem{que}{Problem}

\newtheorem{rem}{Remark}
\newtheorem{lem}{Lemma}
\newtheorem{thm}{Theorem}
\title{Special functions, integral equations and Riemann-Hilbert problem}
\author{  R. Wong$^a$ and Yu-Qiu
Zhao$^b$}

\date{
{\small \emph{
$^a$Department of Mathematics, City University of Hong Kong,
 Kowloon, Hong Kong}}\\
 {\small \emph{$^b$Department of
Mathematics, Sun Yat-sen University, GuangZhou 510275,  China}}}

\begin{document}

\maketitle

   {\bf{Abstract.}}
We consider a pair of   special functions, $u_\beta$ and $v_\beta$,  defined respectively as the   solutions to the integral equations
\begin{equation*}
    u(x)=1+\int^\infty_0 \frac {K(t) u(t) dt}{t+x} ~~\mbox{and}~~v(x)=1-\int^\infty_0 \frac{ K(t) v(t) dt}{t+x},~~x\in [0, \infty),
\end{equation*}
where $K(t)= \frac {1} \pi \exp \left (- t^\beta \sin\frac {\pi\beta} 2\right
)\sin \left ( t^\beta\cos\frac{\pi\beta} 2  \right )$ for  $\beta\in (0, 1)$.
In this note, we establish the existence and  uniqueness   of $u_\beta$ and $v_\beta$ which are bounded and continuous in    $[0, +\infty)$.
Also, we show that a solution to a model Riemann-Hilbert problem   in Kriecherbauer and McLaughlin
[{\it{Int.\;Math.\;Res.\;Not.}},\;1999]
can be constructed explicitly in terms of these functions.
A preliminary asymptotic study is carried out on the Stokes phenomena of these functions  by making use of their  connection formulas.

Several open questions are also proposed for  a  thorough investigation of the analytic and asymptotic properties of the functions $u_\beta$ and $v_\beta$,   and a  related   new special function
 $G_\beta$.



\section{Introduction}\label{introduction}

  Kriecherbauer and McLaughlin \cite{kriecherbauer-mcLaughlin} used the Riemann-Hilbert approach to study the strong asymptotics of the polynomials orthogonal with respect to the
 Freud weight
\begin{equation}\label{freud-weight}w_\beta(x) dx =e^{-\kappa_\beta |x|^\beta} dx,~~x\in \mathbb{R},~~\kappa_\beta =  \frac {\Gamma(\frac \beta 2)\Gamma(\frac 1 2)} {\Gamma(\frac {\beta+1} 2)}, ~~\beta>0.\end{equation}
When $0<\beta <1$, several representative results are obtained  in terms of
a model  Riemann-Hilbert problem (RHP) for a certain $2\times 2$ matrix-valued function $L$:
\begin{align}
   &  L: \mathbb{C}\backslash \mathbb{R}\rightarrow \mathbb{C}^{2\times 2}~ \mbox{is~analytic}, \nonumber\\
  & L_+(s)=L_-(s) v_L(s)~\mbox{for} ~s\in \mathbb{R},\label{L-jump}\\
  & L(s)=O(\ln |s|)~\mbox{as}~s\rightarrow 0,\label{L-behavior-origin}\\
  & L(s)\rightarrow I~\mbox{for}~s\rightarrow \infty,\label{L-behavior-infty}
\end{align}
where the jump matrix is given by
\begin{equation*}
v_L(s)= \left(
                   \begin{array}{cc}
                   1  & -\eta_L(s) \\
                   \eta_L(-s)   & 1\\
                   \end{array}
                 \right),~~s\in  \mathbb{R},
\end{equation*}
with
\begin{equation*}
\eta_L(s)=2i (-1)^{n+1} e^{-|s|^\beta \sin(\pi
\beta/2)}\sin\left(|s|^\beta\cos\frac{\pi\beta}2\right )
{\bf{1}}_{[0, \infty)}(s), ~s\in \mathbb{R},
\end{equation*}and $n$ being the polynomial degree; see  \cite[(6.41)-(6.43)]{kriecherbauer-mcLaughlin}.

Later in this note, we will give   a more natural description of the behavior of $L(s)$ at the origin. Instead of \eqref{L-behavior-origin}, we   henceforth    assume
\begin{equation}\label{L-behavior-origin-new}
L(s)=O\left ( \epsilon_\beta(s)  \right )~~\mbox{as}~~s\to 0,~~\epsilon_\beta(s)= \left\{ \begin{array}{ll}
                                                                                            1, &   1 /2 <\beta<1,\\
                                                                                            \ln |s|, & \beta= 1/ 2, \\
                                                                                           |s|^{\beta-  1 /2}, & 0<\beta< 1/2.
                                                                                          \end{array}
\right.
\end{equation}

{\begin{rem}\label{remark-local-behavior}
The uniqueness of the solution to the RHP for $L$, subject to either   \eqref{L-behavior-origin} or  \eqref{L-behavior-origin-new}, can be justified by using Liouville's theorem. Now let  $L_{\ref{L-behavior-origin}}(s)$ be the solution to the RHP    with local behavior   \eqref{L-behavior-origin}, and $L_{\ref{L-behavior-origin-new}}(s)$ be that with local behavior  \eqref{L-behavior-origin-new}, then
$L_{\ref{L-behavior-origin}}(s) L_{\ref{L-behavior-origin-new}}^{-1} (s)$
has no jump on the real line $\mathbb{R}$, and has only an isolated singularity at the origin. Since the singularity at the origin is not a pole or an essential singularity, it must be removable. In addition, since
$L_{\ref{L-behavior-origin}}(s) L_{\ref{L-behavior-origin-new}}^{-1} (s)$  approaches $I$ as $s\to\infty$,   by Liouville's theorem this matrix-valued   entire function  is the identity matrix.  Thus $L_{\ref{L-behavior-origin}}(s)\equiv L_{\ref{L-behavior-origin-new}}(s)$ if both solutions exist.
In view of Lemma \ref{behavior-at-origin} below, the local behavior  in \eqref{L-behavior-origin-new} is in a sense more natural. Later in Section \ref{sec:boundedness}, we will see from Theorem \ref{boundedness} that we actually have
$L(s)=O(1)$ as $s\to 0$ for the parametrix constructed in \eqref{model-problem-solution}. Accordingly,    both \eqref{L-behavior-origin}  and   \eqref{L-behavior-origin-new} are fulfilled.
\end{rem}
}

In \cite{kriecherbauer-mcLaughlin}, the existence and uniqueness of $L$ were established in   the Wiener class. For $0<\beta <1$,  the leading
 order behavior of the solution to the model problem  {at the origin}   has  not been determined explicitly; see the   comment  in  Deift {\it et al.}\;\cite[p.62, Remark 3.6]{deift-et-al}.
For a brief historical account of the  asymptotics of the orthogonal polynomials associated with  the Freud weight,
the reader is referred to \cite{kriecherbauer-mcLaughlin} and the references therein.

Recently, there are other orthogonal polynomials with weights having  similar singularities as the Freud weight.
For example,
Chen and Ismail \cite{chen-ismail} considered the Freud-like  orthogonal polynomials arising from a
recurrence relation.    The polynomials are   related
to the indeterminate moment problems. The  weight function, supported on $\mathbb{R}$,  has a singularity at $0$ which is   similar to
the Freud weight \eqref{freud-weight}; see also \cite[(1.6)]{dai-ismail-wang} for the weight, and for a  difference equation approach to obtain the Plancherel-Rotach asymptotics.

In \cite{deano-kuijlaars-roman},   Dea\~{n}o,   Kuijlaars  and   Rom\'{a}n  consider the zero distribution of polynomials orthogonal with respect to the Bessel function $J_\nu(x)$. They encountered    a varying exponential weight of Freud type $e^{n\pi|x|}$   with a potential function $\pi |x|$. To justify the existence of a certain local parametrix,
a quite  unnatural restriction on $\nu$ has to be  brought in.  Again, the leading
 order behavior of the parametrix   has   not been   explicitly given.

To solve those parametrices and model problems, one may either  represent them in terms of  known special functions, or  alternatively  use new special functions to serve the same purpose.

The objective of the present note is to introduce a pair of   special functions, $u_\beta$ and $v_\beta$,     determined by  two linear integral equations.
Based on  these functions, the solution to the
 model problem of Kriecherbauer and McLaughlin \cite{kriecherbauer-mcLaughlin}    can readily be constructed.

In the rest of this note, we will   address the unique solvability  of the integral equations in $C[0, \infty)$, derive the behavior of the specific solutions  $u_\beta$ and $v_\beta$ at infinity and at the origin, and carry out  a brief discussion of the Stokes phenomena of these functions, taking the integral equations as resurgent relations.
We will also  propose   a  thorough  investigation of the analytic and asymptotic properties of the new special  functions $u_\beta$ and $v_\beta$,  and a  relevant  new special function
 $G_\beta$.

\section{Integral equations and the RH problem} \label{ie-RHp}

First, let us  derive the integral equations from the model problem for $L$ formulated in Section \ref{introduction}.
From the jump condition \eqref{L-jump} and the behavior \eqref{L-behavior-infty}-\eqref{L-behavior-origin-new}, it is readily seen  that the $(1,1)$-entry $L_{11}(s)$ is analytic in   $ \mathbb{C}\setminus (-\infty, 0]$, such that
\begin{align*}
&\left (L_{11}\right )_+(s)-\left (L_{11}\right )_-(s)=\eta_L(-s)L_{12}(s), ~s\in (-\infty, 0), \\
&L_{11}(s)=1+o(1), ~s\to\infty,  \\
&L_{11}(s)=O(\epsilon_\beta(s)), ~s\to 0.
\end{align*}
Hence, in view of the behavior of $L_{12}$; cf.  \eqref{L-behavior-infty}-\eqref{L-behavior-origin-new},  we have
\begin{equation}\label{L11-integral}
L_{11}(s)=1+\frac 1 {2\pi i}  \int^0_{-\infty} \frac { \eta_L(-\tau)L_{12}(\tau) d\tau}{\tau-s}= 1+\frac 1 {2\pi i}  \int_0^{\infty} \frac { \eta_L(\tau)L_{12}(-\tau) d\tau}{-\tau-s}
\end{equation}for $s\in    \mathbb{C}\setminus (-\infty, 0]$, especially for $s\in (0, \infty)$.
Similarly,   it is also seen from  \eqref{L-jump} that $L_{12}(s)$ is analytic in   $ \mathbb{C}\setminus [0, \infty)$, and solves the scalar RH problem
\begin{align*}
&\left (L_{12}\right )_+(s)-\left (L_{12}\right )_-(s)=-\eta_L(s)L_{11}(s), ~s\in (0, \infty), \\
&L_{12}(s)=o(1), ~s\to\infty,  \\
&L_{12}(s)=O(\epsilon_\beta(s)), ~s\to 0.
\end{align*}
Hence,  we have
\begin{equation} \label{L12-integral}
L_{12}(s)=\frac 1 {2\pi i}  \int_0^{\infty} \frac { -\eta_L(\tau)L_{11}(\tau) d\tau}{\tau-s}~~\mbox{for}~~s\in    \mathbb{C}\setminus [0, \infty),
\end{equation}especially for $s\in (-\infty, 0)$. For $x\in (0, +\infty)$,
now we define
\begin{equation}\label{u-beta-v-beta}
u_\beta(x)=L_{11}(x)+(-1)^n L_{12}(-x)  ,~~v_\beta(x)=L_{11}(x)-(-1)^n L_{12}(-x).
\end{equation}From \eqref{L11-integral}-\eqref{L12-integral} and \eqref{u-beta-v-beta}, it is readily verified that
$u_\beta$ and $v_\beta$ solve, respectively, the following
 integral equations
\begin{align}
&u(x)=1+\int^\infty_0 \frac {K(t) u(t) dt}{t+x},~ ~x\in (0, \infty),\label{u-integral-equation} \\
&v(x)=1-\int^\infty_0 \frac{ K(t) v(t) dt}{t+x},~ ~x\in (0, \infty),\label{v-integral-equation}
\end{align}
where
\begin{align}
          K(t) &=   \frac 1 {2\pi i} \left [
          \exp \left ( e^{(\frac \pi 2 +\frac
{\beta \pi} 2) i} t^\beta\right )
         - \exp \left ( e^{-(\frac \pi 2 +\frac {\beta
\pi} 2) i} t^\beta\right )
\right ] \label{kernel-for-continuation}  \\[.3cm]
           & =    \frac {1} \pi \exp \left (- t^\beta \sin\frac {\pi\beta} 2\right
)\sin \left ( t^\beta\cos\frac{\pi\beta} 2  \right ),\label{kernel-real}
        \end{align}
 and
  $\eta_L(t)=2\pi i (-1)^{n+1} K(t)$ for $t\in [0, +\infty)$.

  It is worth noting that similar coupled scalar integral equations have   been derived in \cite{deano-kuijlaars-roman}, though the equations in that paper were used only to construct a contraction mapping.
Also, it is mentioned in   Fokas {\it et al.}\;\cite[p.161]{fokas-et-al}  that  a   function $u(x)$ can be parametrized  via the solution of a linear integral equation. This is exactly what we are doing: We are defining a pair of new special functions,  using the above integral equations, to construct a parametrix.

 Conversely,  if $u_\beta(x)$ and $v_\beta(x)$ do solve the integral equations, with $u_\beta-1,\; v_\beta-1\in L_2[0, \infty)$, then { we can deduce}
  that $u_\beta(z),\; v_\beta(z)=1+O(1/z)$ as $z\to\infty$ for $|\arg z|<3\pi/2$ and are of the order  $O(\epsilon_\beta(z))$ as $z\to 0$ for $|\arg z|\leq\pi$; see the discussion in the next section, and especially Lemma \ref{behavior-at-origin}. Thus,  we have
 \begin{thm}\label{prop-construction-of-L}The piece-wise analytic function
 \begin{equation}\label{model-problem-solution} L(s)=
\left\{
\begin{array}{lr}
\left(
    \begin{array}{cc}
    L_\beta(s )        &  (-1)^n U_\beta(se^{ - {\pi}   i}) \\
       (-1)^{n} U_\beta(s )  &  L_\beta(se^{-  {\pi}  i}) \\
    \end{array}
  \right), &   0<\arg s < \pi  ;
\\
  &\\
  \left(
    \begin{array}{cc}
     L_\beta(s)    & (-1)^n U_\beta(se^{ \pi   i})\\
(-1)^n U_\beta(s)  & L_\beta(se^{  \pi   i})  \\
    \end{array}
  \right), & -\pi< \arg s < 0  \end{array}
  \right .
\end{equation}
solves the RH problem \eqref{L-jump}, \eqref{L-behavior-infty} and \eqref{L-behavior-origin-new},
where
\begin{equation}\label{L-beta-U-beta}
L_\beta(z)=\frac 1 2 \left ( u_\beta(z)+v_\beta(z)\right ),~~U_\beta(z)=\frac 1 2 \left ( u_\beta(z)-v_\beta(z)\right ), ~~\arg z \in (-\infty, \infty),~z\not=0.
 \end{equation}
 \end{thm}

\section{Unique solvability of the integral equations }

In \cite{kriecherbauer-mcLaughlin}, the unique solvability of the model problem for $L$ is justified in the Wiener class; i.e.,   the Fourier transform of each entry is in $L_1(\mathbb{R})$.

We note that in view of \eqref{L-beta-U-beta} the unique solvability of the integral equations \eqref{u-integral-equation}-\eqref{v-integral-equation}
can accordingly be deduced from that of $L$.

Now we take an alternative approach, using the integral equations alone. We first show that both $u_\beta(x)-1$ and $v_\beta(x)-1$ are in    $L_2[0, +\infty)$. We begin with the following lemma:
\begin{lem}\label{compact-operator}
   The operator $T$ is a compact operator on $L_2[0, +\infty)$, where
   \begin{equation}\label{compact-T}
    (Tu)(x)=\int^\infty_0\frac {K(t)u(t) dt}{t+x}.
   \end{equation}
\end{lem}\vskip .2cm

\noindent
{\bf{Proof. }} The lemma is easily justified by noting that $\frac {K(t)}{t+x}\in L_2{([0,+\infty)\times [0,+\infty))}$. Thus, $T$ is a Hilbert-Schmidt operator, and hence  is compact; cf. e.g. \cite[p.277]{yosida}.\hfil\qed\vskip.4cm

The next step is to show that the operator $I-T$ has trivial null space. To this aim, we need to know more about the analytic structure of the solutions to \eqref{u-integral-equation} and \eqref{v-integral-equation}. Indeed, if $u(x)$ solves \eqref{u-integral-equation}, then we see from   equation \eqref{u-integral-equation}  that $u$ can be extended analytically to $\mathbb{C}\setminus (-\infty, 0]$, and that
$u(z)=1+O(1/z)$ in subsectors. From \eqref{kernel-for-continuation}, we see that the kernel $K(t)$ can be extended analytically to  complex $t$ with $\arg t\in (-\infty, \infty)$, and is exponentially small for $|\arg t |<\pi/2$.
Using   Cauchy's theorem, it can be shown that
\begin{equation}\label{analytic-continuation-rotation}
u(z)= 1+ \int^\infty_0 \frac {K(te^{i\theta}) u(te^{i\theta}) dt}{t+e^{-i\theta} z},~~|\theta|<\frac \pi 2,
\end{equation}
which enables us to extend $u(z)$ to the larger sector $|\arg z|<3\pi/2$, and
\begin{equation}\label{u-behavior-at-infty}
u(z)=1+O(1/z),~~\mbox{as}~z\to\infty,~~|\arg z|<3\pi/2.
\end{equation}Similar results can be obtained for $v(z)$.
This   procedure   can be continued,  since by using the Plemelj formula we have from \eqref{u-integral-equation}-\eqref{v-integral-equation}
\begin{equation}\label{connection-u-v} u(ze^{-\pi i})-u(ze^{\pi i})= 2\pi i K(z) u(z),~~v(ze^{-\pi i})-v(ze^{\pi i})= -2\pi i K(z) v(z),\end{equation}
initially for real and positive $z$, and then for complex  $z$.
Appealing to the formulas in \eqref{connection-u-v}, we see that the above process of analytic  continuation  can keep on  going until we have the results valid for    $\arg z\in \mathbb{R}$.
For later use, we need the  following  behavior of $(Tu)(z)$ at the origin:
\begin{lem}\label{behavior-at-origin}
   For $u\in L_2[0, +\infty)$, we have
    $(Tu)(z)=O(\epsilon_\beta(z))$ for $|\arg z |\leq \pi$  as $z\to 0$; cf. \eqref{L-behavior-origin-new},
  where
   $
    \epsilon_\beta(z)=
     \left\{ \begin{array}{ll}
                                                                                         1, &  1 /2<\beta<1,\\
                                                                                            \ln |z| , & \beta= 1/ 2, \\
                                                                                            |z|^{\beta- 1/ 2} , & 0<\beta<  1 /2.
                                                                                          \end{array}
\right.
    $
\end{lem}\vskip .2cm

\noindent
{\bf{Proof. }}
    Recall the function $K(t)$ in \eqref{kernel-real}. For real and positive $z$, we have  by using H\"{o}lder's inequality and splitting the   interval of integration
\begin{equation*}
  | (Tu)(z)|^2\leq \|u\|_2^2 \int^\infty_0 \frac {K(t)^2 dt}{(t+z)^2} = O\left ( \int^1_0 \frac {t^{2\beta} dt}{ (t+z)^2}\right ) +O(1) = O(\epsilon_\beta(z))
\end{equation*}as $z\to 0$. The last equality can be obtained by calculating   a  residue.   With slight modifications, and taking into account the rotation formula \eqref{analytic-continuation-rotation}, we can extend the above inequality to  $|\arg z |\leq \pi$.\hfill\qed\vskip.4cm

We show  that the operator $I-T$ has trivial null space,  by using   a vanishing lemma technique.  Similar argument can be found in, e.g., \cite{kriecherbauer-mcLaughlin}.
\begin{lem}\label{vanishing-lemma}
 Assume that   $u\in  L_2[0,+\infty)$ and  $u-Tu=0$ in  $L_2[0,+\infty)$. Then it holds  $u\equiv 0$ for $x\in [0, \infty)$ in  $L_2[0,+\infty)$.
\end{lem}\vskip .2cm

\noindent
{\bf{Proof. }}
Indeed, from  $u=Tu$ we see that the above analytic continuation procedure  works also  for the present $u(z)$, the connection formula \eqref{connection-u-v} still holds, and instead of \eqref{u-behavior-at-infty} we have
$u(z)=O(1/z)$ as $z\to\infty$ for $|\arg z|\leq \pi$. We introduce an auxiliary vector  function
\begin{equation*}
  A(z)=\left\{ \begin{array}{ll}
               \left(
         \begin{array}{cc}
           u(z), & u(ze^{-\pi i})
         \end{array}
       \right),  & \arg z\in (0, \pi), \\[.2cm]
                 \left(
                   \begin{array}{cc}
                     u(z), & u(ze^{\pi i})
                   \end{array}
                 \right),
                  & \arg z\in (-\pi, 0).
               \end{array}\right .
\end{equation*}
It is readily seen that $A(z)$ is analytic in both the upper and lower half plane,  behaves uniformly as $O(1/z)$ at infinity, and has jumps along the real axis as
\begin{equation*}
  A_+(x)=A_-(x)J_A(x),~~J_A(x)=\left\{\begin{array}{ll}
                      \left(
                  \begin{array}{cc}
                    1 & 2\pi i K(x) \\
                    0 & 1 \\
                  \end{array}
                \right)     &\mbox{for}~ x\in (0, +\infty), \\[.4cm]
                    \left(
                  \begin{array}{cc}
                    1 & 0\\
                    -2\pi i K(|x|) & 1 \\
                  \end{array}
                \right)       &  \mbox{for}~ x\in (-\infty, 0);
                      \end{array}\right .
\end{equation*}cf. \eqref{connection-u-v}.

Denote by $A^*$ the   conjugate transpose of $A$. Since  $A_+A_-^*$ has an analytic continuation to the upper half plane, and
$A_-A_+^*$ to the lower half, in view of the behavior $A(z)=O(1/z)$ at infinity, by the Cauchy integral theorem  we have
 \begin{align*}
           &\int_{\mathbb{R}} A_-(x)J_A(x) A_-^*(x) dx  =  \int_{\mathbb{R}} A_+(x)  A_-^*(x) dx=0,\\[.2cm]
&\int_{\mathbb{R}} A_-(x)J_A^*(x) A_-^*(x) dx  =  \int_{\mathbb{R}} A_-(x)  A_+^*(x) dx=0.
        \end{align*}
Here,  use has also been made of the fact that $u(z)=O  (\epsilon_\beta(z))$ for small $z$ with  $|\arg z |\leq \pi$; see Lemma \ref{behavior-at-origin}.
Summing up the integrals gives
\begin{equation}\label{integral-equal-to-zero}
\int_{\mathbb{R}} A_-(x)\left (J_A(x)+ J_A^*(x)\right ) A_-^*(x) dx   =0.
\end{equation} It can be seen from \eqref{kernel-real} that $|\pi i K(|x|)|<1$ for $x\in \mathbb{R}$. Hence,
  \begin{equation*}
    J_A(x)+ J_A^*(x)=\left(
                  \begin{array}{cc}
                    2 & 2\pi i K(|x|)\\
                    -2\pi i K(|x|) & 2 \\
                  \end{array}
                \right)
  \end{equation*}
   is positive definite  for $x\in  \mathbb{R}$. From \eqref{integral-equal-to-zero} one deduces
that $A_-(x)\equiv 0$, and thus $u(x)\equiv 0$ for $x\in (0, +\infty)$.\hfill\qed\vskip .4cm

We note   that $\int^\infty_0\frac {K(t) dt}{t+x}\in L_2[0, +\infty)$ since it is smooth on $(0, \infty)$, bounded at $0^+$, and of order $O(1/x)$ at $+\infty$. Hence,      a combination  of Lemmas \ref{compact-operator} and \ref{vanishing-lemma}  gives  the unique solvability of the following integral equation
\begin{equation*}
 u_1(x)-Tu_1(x) =\int^\infty_0 \frac {K(t) dt}{t+x}
\end{equation*}in the space  $L_2[0, +\infty)$.   Similar analysis can be carried out for the other equation \eqref{v-integral-equation}. To conclude, we have the following:
 \begin{thm}\label{unique-solvability}  There exist  unique solutions $u(x)$ and $v(x)$  to the integral equations \eqref{u-integral-equation} and \eqref{v-integral-equation}, respectively,  such that $u(x)-1\in L_2[0, +\infty)$ and $v(x)-1\in L_2[0, +\infty)$.\hfill\qed
\end{thm}

\section{Boundedness of $u_\beta(x)$ and $v_\beta(x)$ on $[0, \infty)$}\label{sec:boundedness}
The main result in the last section is that $u_\beta(x)-1\in L_2[0, +\infty)$. By using successive approximation, one can actually show that $u_\beta(x)$ is continuous and bounded for $\beta> {\beta_c:=} 0.4158853544$.
{Indeed, it suffices to show that the operator $T$ defined in \eqref{compact-T} is a contraction  on the space  $L_\infty$  for those $\beta$. Straightforward estimation gives
\begin{equation*}
   |(Tu)(x)|\leq M_\beta \max_{[0, +\infty)} |u(t)|~\mbox{for}~x\in [0, +\infty),~~\mbox{with}~M_\beta=\int^\infty_0 \frac {|K(t)|dt} t\leq \frac {\cot(\frac {\pi\beta} 2)}{\pi\beta}.
\end{equation*}
 The bound  ${\cot (\frac {\pi\beta} 2)}/{\pi\beta}$  is monotonically  decreasing for $\beta\in (0,1)$, and  is less than  $1$ for $\beta>\beta_c$. }

In  this section, we shall show that {$u_\beta(x)$ and $v_\beta(x)$ are in fact  bounded}
 for all $\beta\in (0, 1)$.

  \begin{thm}\label{boundedness}  Let $u_\beta(x)$ and $v_\beta(x)$ be  solutions of the integral equations in \eqref{u-integral-equation} and \eqref{v-integral-equation}, respectively.
  Then, $u_\beta(x),~v_\beta(x)\in L_\infty[0, \infty)$ and $u_\beta(x), ~v_\beta(x)\in C[0, \infty)$.
  \hfill\qed
\end{thm}

\noindent
{\bf{Proof. }} Let $\tilde u(x)=u_\beta(x)-1$. By Theorem \ref{unique-solvability}, $\tilde u(x)\in L_2[0, \infty)$. Rewrite equation \eqref{u-integral-equation} as
\begin{equation*}
    \tilde u(x)= \int^\infty_0 \frac {K(t) \tilde u(t) dt}{t+x}+u_0(x),~~x\in [0, \infty),
\end{equation*}where
\begin{equation*}
u_0(x)=\int^\infty_0 \frac {K(t)}{t+x}  dt.
\end{equation*}Clearly,
$u_0(x) \in  L_2[0, \infty)\cap L_\infty[0, \infty)$. The function  $\tilde u(x)$ is smooth in $(0, \infty)$, since it is an analytic function in the cut plane $\mathbb{C}\setminus (-\infty, 0]$.
To prove that   $\tilde u(x)$ is also   bounded there, we divide our discussion into several  cases.
In general, for all $\beta\in (0, 1)$, we have
\begin{equation*}
  |\tilde u(x)|\leq |u_0(x)|+ \int^\infty_0 \frac {|K(t)| \; |\tilde u(t)| dt}{t+x}\leq  \frac 1 x \int^\infty_0  |K(t)| (1+ |\tilde u(t)|) dt \leq   \frac 1 x\left ( {\|K\|_1+\| K\|_2 \|\tilde u\|_2} \right )
\end{equation*}for all positive $x$, and especially for $x\geq 1$. We may focus on $0<x\leq 1$.
\vskip .5cm

{\bf{Case}} (i)  $\frac 1 2<\beta <1$. For $0<x\leq 1$, we have
\begin{equation*}
  |\tilde u(x)|\leq a+ \int^\infty_0 \frac {|K(t)| \; |\tilde u(t)|}{t+x} dt\leq  a+ \int^\infty_0  \frac {|K(t)|} t  |\tilde u(t)| dt\leq  a+ \left \|\frac { K(t) } t\right\|_2 \|\tilde u\|_2,
\end{equation*}   where $a=\|u_0\|_\infty$. Thus, $\tilde u(x)$ is bounded. Here, we have used the fact that  $\frac { K(t) } t\in L_2[0, \infty)$.
\vskip .5cm

{\bf{Case}} (ii)   $\frac 1 4<\beta \leq \frac 1 2$.  In this case, we   write $\beta=\frac 1 4+\varepsilon$ with $\varepsilon\in (0, \frac 1 4]$.    We
first recall the inequality \cite[p.30]{mitrinovic}
\begin{equation}\label{A}
 x^{\frac 1 p}y^{\frac 1 q} \leq \frac x p+\frac y q\leq x+y,
\end{equation}where $x\geq 0$, $y\geq 0$ and $1/p+1/q=1$ with $p>1$. From this inequality, we have
\begin{equation}\label{B}
 t+x\geq t^{\beta+\frac 1 2-\frac \varepsilon 2} x^{\frac 1 2-\beta+\frac \varepsilon 2} .
\end{equation}Note that
  $\beta+\frac 1 2 -\frac \varepsilon 2=\frac 3 4+\frac  \varepsilon 2>0$  and $\frac 1 2-\beta + \frac  \varepsilon 2\geq \frac  \varepsilon 2>0$. Consequently, we have
\begin{equation}\label{estimate}
  |\tilde u(x)|\leq a+ \int^\infty_0 \frac {|K(t)| \; |\tilde u(t)|}{t+x} dt\leq    a+ \left \|\frac { K(t) } {t^{\beta+\frac 1 2-\frac \varepsilon 2}}\right\|_2 \|\tilde u\|_2\frac 1 {x^{\frac 1 2-\beta+\frac \varepsilon 2}}\leq \frac c {x^{\frac 1 2-\beta+\frac \varepsilon 2}},~~x\in (0, 1],
\end{equation}where $c$ is a generic positive constant.
The fact that ${ K(t) } /{t^{\beta+\frac 1 2-\frac \varepsilon 2}}\in L_2[0, \infty)$ is readily justified by its behavior $O(t^{\frac \varepsilon 2-\frac 1 2})$ at $t=0^+$, and its exponential decay at $t=+\infty$.
This estimate can be used to  self-improve and  achieve
\begin{equation*}
  |\tilde u(x)|\leq a+ \int^\infty_0 \frac {|K(t)| \; |\tilde u(t)|}{t}dt\leq  a+ \left \|\frac { K(t) } {t^{1-\beta+  \varepsilon  }}\right\|_2 \|t^{-\beta+\varepsilon} \tilde u(t) \|_2
\end{equation*}   for  $x\in (0, \infty)$.  To see that the two quantities in the last term of the above inequality are finite, we just check their behavior as
$t\to 0^+$ and  as $t\to +\infty$. For instance,  as $t\to 0^+$,  we have $  { K(t) } /{t^{1-\beta+  \varepsilon  }}= O  ( t^{ 2\beta-1-\varepsilon}   )=O  ( {t^{\varepsilon-\frac 1 2 }}  )$ since $\beta=\frac 1 4+\varepsilon$ and
$t^{-\beta+\varepsilon} \tilde u(t)= O (  {t^{\frac \varepsilon 2-\frac 1 2  }}  )$ by the estimate in \eqref{estimate}.
\vskip .5cm

{\bf{Case}} (iii)   $\frac 1{2(k+1)} <\beta \leq \frac 1 {2k}$, $k=2,3,\cdots$.  The analysis  is very much the same as that given  for  the previous case, except that  we need to repeat   the self-improvement argument several  times. First,
let us write $\beta= \frac 1{2(k+1)}+\varepsilon$ with $\varepsilon\in (0, \frac 1 {2k(k+1)}]$.
Straightforward calculation shows
\begin{equation*}
k\beta+\frac 1 2-\frac {(2k-1)\varepsilon }2 > \cdots > 2\beta+\frac 1 2-\frac {3\varepsilon }2 > \beta+\frac 1 2-\frac { \varepsilon }2 =\frac 1 {2(k+1)}+\frac 1 2+\frac \varepsilon 2>0,
\end{equation*}
and
 \begin{equation*}
\frac 1 2- \beta+ \frac { \varepsilon }2 > \frac 1 2- 2\beta+ \frac {3\varepsilon }2>\cdots >
\frac 1 2 -k\beta+ \frac {(2k-1)\varepsilon }2  =\frac 1 {2(k+1)}-\frac \varepsilon 2\geq \frac {2k-1}{4k(k+1)} >0.
\end{equation*}
Hence, for $x\in (0, 1]$, repeated application  of  the inequality in \eqref{A} and the self-improvement steps in \eqref{B} and \eqref{estimate} lead to
\begin{equation*}\left\{
 \begin{array}{l}
     t+x \geq t^{\beta+\frac 1 2-\frac \varepsilon  2} x^{\frac 1 2-\beta+\frac \varepsilon  2}~~~\mbox{and}~~~   \\
     |\tilde u(x)|\leq a +  { \left \| \frac {K(t)} {t^{\beta+\frac 1 2-\frac \varepsilon 2}}\right\|_2 \|\tilde u(t)\|_2}\;  \frac 1 {x^{\frac 1 2-\beta+\frac \varepsilon  2}}\leq
 \frac c {x^{\frac 1 2-\beta+\frac \varepsilon  2}},
\end{array}
\right .
\end{equation*}
 \begin{equation*}\left\{ \begin{array}{l}
                             t+x \geq t^{2\beta+\frac 1 2-\frac {3\varepsilon }2} x^{\frac 1 2-2\beta+\frac {3\varepsilon}2}~~~\mbox{and}\\
                             |\tilde u(x)|
\leq a +   { \left \| \frac {K(t)} {t^{\beta+\frac 1 2-\frac \varepsilon 2}}\right\|_2 \left \|t^{-\beta+\varepsilon} \tilde u(t)\right \|_2}\; \frac 1 {x^{\frac 1 2-2\beta+\frac {3\varepsilon } 2}}
 \leq \frac c {x^{\frac 1 2-2\beta+\frac {3\varepsilon }2}},
                          \end{array}
    \right .
 \end{equation*}
   \begin{equation*}
             ~~~~ \vdots
       \end{equation*}
\begin{equation*}\left\{
 \begin{array}{l}
   t+x \geq t^{k\beta+\frac 1 2-\frac {(2k-1)\varepsilon }2} x^{\frac 1 2-k\beta+\frac {(2k-1)\varepsilon }2}~~~\mbox{and} \\
   |\tilde u(x)| \leq a +  { \left \| \frac {K(t)} {t^{\beta+\frac 1 2-\frac \varepsilon 2}}\right\|_2 \left \|t^{-(k-1)(\beta-\varepsilon)} \tilde u(t)\right \|_2}\;  \frac 1 {x^{\frac 1 2-k\beta+\frac {(2k-1)\varepsilon }2}}   \leq \frac c {x^{\frac 1 2-k\beta+\frac {(2k-1)\varepsilon }2}}.
 \end{array}
\right .
\end{equation*}
The last estimate ensures the finiteness of the two $L^2$-norms in the following inequalities
\begin{equation*}
  |\tilde u(x)|\leq a+ \int^\infty_0 \frac {|K(t)| \; |\tilde u(t)|}{t} dt\leq a+ \left \|\frac { K(t) } {t^{1-k\beta+   k \varepsilon    }}\right\|_2 \|t^{-k\beta+ k\varepsilon  } \tilde u(t) \|_2 .
\end{equation*}  This establishes the boundedness of $u_\beta (x)$ for  $x\in (0, 1]$, and hence also  for $x\in (0, \infty)$.
Since $u_\beta(x)$ is bounded on $(0, \infty)$, the continuity of $u_\beta(x)$ at $x=0$ from the right now follows from equation \eqref{u-integral-equation} and the Lebesgue dominated convergence theorem.

An analogous argument gives $v_\beta(x)\in L_\infty [0, \infty)$ and  $v_\beta(x)\in C [0, \infty)$. \hfil\qed

\section{Asymptotic behavior and Stokes phenomenon}

From the integral equations,  and in view of \eqref{analytic-continuation-rotation}, we readily obtain the asymptotic approximations
\begin{equation}\label{u-beta-v-beta-infty}
 u_\beta (z)\sim 1+\sum^\infty_{k=1} \frac {c_k}{z^k} ~~\mbox{and}~~ v_\beta(z)\sim 1+\sum^\infty_{k=1} \frac {d_k}{z^k}~~\mbox{as}~|z|\to+\infty,~~|\arg z| <   \frac {3\pi} 2,
\end{equation}
where
\begin{equation*}
  c_k=(-1)^{k-1} \int ^\infty_0 K(t) t^{k-1}u_\beta(t)dt~~\mbox{and}~~d_k=(-1)^{k} \int ^\infty_0 K(t) t^{k-1}v_\beta(t)dt
\end{equation*}for $k=1,2,\cdots$;
cf.  Boyd \cite{boyd} for asymptotic approximations of Stieltjes transforms. See also Wong \cite{wong-book}, and Wong and Zhao \cite{wong-zhao_S-T}.

To satisfy certain    jump conditions and   matching conditions,
the Stokes phenomenon plays a crucial role in the construction of local parametrices for RHPs.
Now we turn to a brief discussion  of the  Stokes phenomena of $u_\beta$ and $v_\beta$, based on the connection formulas in \eqref{connection-u-v}, and regarding    equations \eqref{u-integral-equation}-\eqref{v-integral-equation} as   exact resurgence relations. A similar analysis  was carried out in Boyd \cite{boyd} for the modified Bessel function; see also Xu and Zhao \cite{xu-zhao_Resurgence} for an attempt to link together the RH approach with   resurgence properties.

Assume that  $u(x)=u_\beta(x)$ solves the integral equation \eqref{u-integral-equation}  with  $u-1\in L_2[0, \infty)$, and $u(z)$ denotes its analytic continuation. In general, we see that $u(z)= 1+O(1/z)$ for large $z$ as $|\arg z|<\frac {3\pi}2$; cf. \eqref{u-behavior-at-infty}.
 Coupling  \eqref{u-behavior-at-infty} and \eqref{connection-u-v}, one can verify that
 \begin{equation*}
   u(z)=1+O\left (\frac 1 z\right ) +  \exp\left ( e^{-i  \pi (1+ \beta )/  2} (z e^{-\pi i})^\beta\right )
 \end{equation*} for $\frac \pi 2<\arg z\leq \frac {3\pi} 2$; see also \eqref{kernel-for-continuation}. The term on the extreme right is exponentially large, when $\arg z$ goes beyond $\frac {3\pi} 2$.  According to   Berry's Stokes smoothing, the switch-on of the exponential occurs when   $1+O(1 /z)$ is the most dominant as compared with the exponential term, namely,
 when $-\left [\frac \pi 2+\left ( \frac {3\pi }2 -\arg z\right )\beta\right ]=-\pi$, i.e.,       $\arg z =\frac {3\pi} 2-\frac \pi {2\beta}$. We give some details in what follows.

In view of \eqref{kernel-for-continuation}, we may write
\begin{equation}\label{u-divide-kernels}
  u(z)=1+ \frac  1 {2\pi i} \int^\infty_0\frac {e^{ i e^{{\beta\pi i} /2  } t^\beta } u(t)    dt}{t +z}-
  \frac 1 {2\pi i} \int^\infty_0\frac {e^{ -i e^{-   {\beta\pi i} /2 } t^\beta  } u(t)    dt}{t +z}:=1+I_P-I_N,
  \end{equation}initially for real positive $z$, and then extended elsewhere.  Now we assume that $\beta\in [1/3, 1)$, and  take $I_N(z)$ as an example. Rotating the integration path clockwise to $\arg t=\phi \in [-\pi, 0)$, and making a change of variables
$t=e^{-\pi i} z \tau^\alpha$    with $\alpha=1/\beta$,   we have
 \begin{equation}\label{stokes-nonstokes}
I_N(z)= \frac \alpha {2\pi i} \int^\infty_0 \frac {e^{-\rho e^{i\theta} \tau} u(e^{-\pi i}z \tau^\alpha) \tau^{\alpha-1} d\tau}{\tau^\alpha-1},
 \end{equation}where $\theta=\beta\left ( \arg z-\frac \pi 2(3-\alpha)\right )$,   the large parameter is  $\rho=|z|^\beta$,     and  the path of integration is indented to pass above the pole at $\tau =1$.
Before reaching the equality in \eqref{stokes-nonstokes}, we have already rotated the path $\arg\tau=\phi+\pi$ to the positive real line.

Substitute  the identity $\frac 1 {\tau ^\alpha-1}=-\sum^{N-1}_{k=1}\tau ^{k\alpha} +\frac {\tau ^{N\alpha}} {\tau ^\alpha-1}$ into the  integral in \eqref{stokes-nonstokes}.
 The remainder in the resulting expansion is an integral with a saddle point as well as a pole. Optimal truncation of the expansion occurs when the saddle point coalesces with the pole (see \cite{boyd,wong-zhao_M-L}); that is, when  $\rho\approx N\alpha$  or, more precisely, $\alpha N=\rho+r$, with $r$ being bounded. The result we have now  is an  exponentially improved asymptotic expansion, whose  error term is given by
\begin{equation*}
   \frac 1 {2\pi i} \int_0^\infty \left \{ \frac {\tau -1}{\tau ^\alpha-1} u(e^{-\pi i}z\tau ^\alpha  ) \alpha \tau ^{\alpha-1+r}  \right \}   \frac {e^{-\rho\left ( e^{i\theta} \tau -\ln \tau \right )} d\tau }{\tau -1},
\end{equation*}which  exhibits a Stokes smoothing  of  error-function type  as $\theta\to 0^-$ when the saddle point $\tau=e^{-i\theta}$ coalesces with the simple pole $\tau=1$.
 This indicates that in our  case a Stokes line  occurs when we have  $\theta=0$, i.e., $\arg z= \frac \pi 2(3-\alpha)$. The analysis for the other integral $I_P(z)$   gives    another  Stokes line $\arg z=-\frac \pi 2(3-\alpha)$.

A slight modification is   needed for the case $\beta\in [1/5, 1/3)$, since $\phi_0=  -\frac \pi 2  (\alpha-1)\in [-2\pi, -\pi)$ in this case. Rotating  the integration path clockwise and picking up the residue $r_N(z)=-\exp\left ( e^{-(\frac \pi 2+\frac {\beta\pi} 2) i} (ze^{-\pi i})^\beta\right )u(ze^{-\pi i})$ at $t= ze^{-\pi i}$, one has
 \begin{equation*}
   I_N(z)=r_N + \frac 1 {2\pi i} \int^{\infty e^{i\phi}} _0\frac {e^{ -i e^{- \frac {\beta\pi i} 2 } t^\beta  } u(t)    dt}{t +z}=J_N(z)+ \frac 1 {2\pi i} \int^{\infty e^{i\phi}} _0\frac {e^{ -i e^{- \frac {\beta\pi i} 2 } t^\beta  } u(te^{2\pi i})    dt}{t +z}   ,
 \end{equation*}
where $\phi<-\pi$,  and $J_N(z)$ is a sum of exponentials and integrals. The second equality is obtained by substituting in $u(t)=u(te^{2\pi i})+K(te^{\pi i})u(te^{\pi i})$; cf. \eqref{connection-u-v}.  Now we focus on the last integral. Making a change of variables $t=e^{-\pi i} z \tau^\alpha$, and using the same optimal truncation $\alpha N=\rho+r$ with $\rho=|z|^\beta$, we obtain  a slightly different  error term
\begin{equation*}
   \frac 1 {2\pi i} \int_0^\infty \left \{ \frac {\tau -1}{\tau ^\alpha-1} u(e^{\pi i}z\tau ^\alpha  ) \alpha \tau ^{\alpha-1+r}  \right \}   \frac {e^{-\rho\left ( e^{i\theta} \tau -\ln \tau \right )} d\tau }{\tau -1},
\end{equation*}
where again $\theta=\beta\left ( \arg z-\frac \pi 2(3-\alpha)\right )$,   the large parameter $\rho=|z|^\beta$,     but  the path of integration is indented to pass below the pole at $\tau =1$.
As before, a   Stokes line in this case is identified as $\arg z=-\frac\pi 2(\alpha -3) \in (-\pi, 0)$.  Another symmetric Stokes line  $\arg z=\frac\pi 2(\alpha-3)$  can be found by analyzing $I_P(z)$; see \eqref{u-divide-kernels}.

In general, for $  \beta\in [\frac 1 {4l+3}, \frac 1 {4l-1})$, $l=1,2,\cdots$,  repeated use of the connection formula \eqref{connection-u-v} gives
 \begin{equation*}
   I_N(z)=J_N(z)+ \frac 1 {2\pi i} \int^{\infty e^{i\phi}} _0\frac {e^{ -i e^{- \frac {\beta\pi i} 2 } t^\beta  } u(te^{2l \pi i})    dt}{t +z}   ,
 \end{equation*}
where $\phi< -(2l-1)\pi$. Here again   $J_N(z)$ is a combination of exponentials and integrals. Similar analysis can be carried out for the last integral to locate the Stokes line.
The same treatment also works  for $v_\beta$. The reader is referred to \cite[p.360]{wong-zhao_M-L}  for full details.


To summarize, we have
\begin{thm}Assume that the functions $u_\beta(z)$ and $v_\beta(z)$, solving respectively \eqref{u-integral-equation} and \eqref{v-integral-equation}, such that $u_\beta-1,~v_\beta-1\in   L_2[0, \infty)$. Then the Stokes lines for both functions are
\begin{equation*}
  \arg z=\pm \left\{  \frac \pi 2(3-\alpha)+2(l-1)\pi\right\},
\end{equation*}
  where  $\alpha=\frac 1 \beta \in ( 4l-3, 4l+1]$, $l=1,2,\cdots$.\hfill\qed
\end{thm}

\section{Discussion}

In the previous sections, we have determined   a pair of
new special functions $u_\beta$ and $v_\beta$ via  integral equations. Preliminary analysis reveals some of the   analytic structures of the functions,  such as their analytic continuations,
   connection formulas, and asymptotic results involving   Stokes phenomenon.  Still we have the following open problems.
\begin{que}
 Here we propose a thorough  investigation of the analytic and asymptotic  properties of the  functions  $u_\beta$ and $v_\beta$,   such as
  their  zeros,  modulus function  \cite[\S 10.18]{nist}, kernel  \cite{kuijlaars-vanlessen}, and differential or difference  equations, etc. For example, one may   begin with
    \begin{description}
    \item[(i)] Determine   the coefficients $c_k$ and $d_k$ of the expansion at infinity given in \eqref{u-beta-v-beta-infty}, which  are now expressed in terms of $u_\beta$ and $v_\beta$. It would be interesting and challenging to decode the coefficients from \eqref{u-integral-equation} and \eqref{v-integral-equation} in an explicit manner;
    \item[(ii)] On account of \eqref{G-origin-rearrange}, one may { reasonably expect} that the behaviors of $u_\beta$ and $v_\beta$ at the origin are of the form $\sum_{k=0}\tilde c_k z^k+\sum_{k=1} \tilde d_k z^{\beta k}$; see also Theorem \ref{boundedness} and Remark \ref{remark-local-behavior} for the boundedness of these functions. A natural problem is to determine explicitly the coefficients  $\tilde c_k$ and  $\tilde d_k$, at least for $k=0$ and $1$.
  \end{description}
\end{que}

The Stieltjes transform of the function $K(t)$ given in \eqref{kernel-for-continuation} suggests that it would  also be  of interest to study  the related integral
\begin{equation}\label{G-definition}
 G_\beta (z):=\int^\infty_0\frac {e^{-t^\beta}}{t+z} dt
\end{equation}for  $\beta >0$,  defined initially for $\arg z\in  (-\pi, \pi)$  and then analytically continued to $\arg z\in  \mathbb{R}$.
The  special case when $\beta=2$ is termed  the Goodwin-Staton integral \cite[(7.2.12)]{nist}; see also Jones \cite{jones}, where, by making use of the incomplete gamma function,
    asymptotic formulas are obtained for a slightly more general    integral.

The Mittag-Leffler function can be represented  in terms of this function.
For $\alpha=1/\beta>0$, the Mittag-Leffler function is given by
\begin{equation}\label{mittag-leffler-path-deformed}
E_\alpha (z):=\sum ^\infty _{n=0} \frac {z^n }{\Gamma(\alpha n+1)}
=\beta \sum e^{Z_k}+     \frac 1 {2\pi i}   \int_C  \frac { t^{\alpha-1}e^t dt } {t^\alpha-z}:=\Sigma(z) +I(z),\end{equation}
where $Z_k=z^\beta e^{2\pi i k \beta}$ for integer $k$  and the summation is taken over all  $k$  satisfying  $| \arg Z_k | <\pi$. The path $C$ is a Hankel contour, and  dented properly  if $ | \arg Z_k | =\pi$ for a certain integer $k$; cf. \cite{wong-zhao_M-L}.

As in  Wong and Zhao \cite{wong-zhao_M-L}, it is appropriate to divide our discussion into several cases: (i)  $\alpha$ is a positive integer,  (ii) $\alpha\in (2l-1, 2l)$, $l=1, 2, 3,\cdots$, and (iii)  $\alpha\in (2l, 2l+1)$, $l=0, 1, 2, \cdots$.

In case (i), i.e., when $\alpha$ is a  positive integer, the contour integral in \eqref{mittag-leffler-path-deformed} vanishes;    thus the function $E_\alpha(z)$ is an elementary function. In   case (ii), we have
\begin{equation*}
I(z)
=    \frac \beta  {2\pi i}  G_\beta ( e^{(2l-1-\alpha)\pi i } z) - \frac \beta  {2\pi i}  G_\beta ( e^{(\alpha -2l+1)\pi i } z)    \end{equation*}
initially for $\arg z=0$, and then extend the region of validity by analytic continuation.  Similarly,  in case (iii),
\begin{equation*}
I(z)
=    \frac \beta  {2\pi i}  G_\beta ( e^{(2l+1-\alpha)\pi i } z) - \frac \beta  {2\pi i}  G_\beta ( e^{(\alpha -2l-1)\pi i } z). \end{equation*}

Asymptotic expansions  of  $G_\beta$ can be derived as follows.
For example,  at infinity,  we have
\begin{equation}\label{G-infinity}
G_\beta(z)\sim \sum_{n=1}^\infty  \frac {(-1)^{n-1} \Gamma(n/\beta)}{\beta z^n}
\end{equation}initially for $|\arg z|<\pi$. The region of validity can be extended    to $|\arg z|<\pi +\frac { \pi} {2\beta}$.
Also, at the origin, we write  $G_\beta(z)=\left\{\int_0^1+\int_1^\infty\right \} \frac {e^{-t^\beta} dt}{t+z}$.  Expanding the exponential in the first integral in convergent series, and using Cauchy residue theorem, we have
 \begin{equation}\label{G-origin-rearrange}
G_\beta(z)\sim -\ln z+
\sum_{n=1}^\infty  \frac {(-1)^{n+1} \pi } {n! \sin(\beta n\pi)} z^{\beta n} +\sum^\infty_{k=0}c_{\beta, k}  (-1)^k  z^k
,~~z\to 0
\end{equation}so long as  $\beta$ is not a rational number,
where
\begin{equation}\label{c-beta-0}
 c_{\beta, 0}=\sum_{n=1}^\infty  \frac {(-1)^n} {n!  n\beta }+   \frac {\Gamma(0, 1)} \beta=\frac 1 \beta \left [ \int^1_0\frac {e^{-t}-1} t dt+\int^\infty_1\frac  {e^{-t}} t dt\right ]=-\frac \gamma \beta
\end{equation}
and
\begin{equation*}
 c_{\beta, k}=
\sum_{n=0}^\infty  \frac {(-1)^n} {n! (n\beta-k) }+   \frac {\Gamma(-  k/ \beta, 1)} \beta,~~k=1,2,\cdots,
\end{equation*}where
 $\Gamma(a, z)=\int^\infty_z t^{a-1} e^{-t} dt$ is the incomplete Gamma function.
For the  last equality in \eqref{c-beta-0},
see  \cite[(6.2.1), (6.2.3)-(6.2.4)]{nist},   $\gamma$ being  Euler's constant.
It can be shown that the expansion in \eqref{G-origin-rearrange} is valid in the sector  $|\arg z|<\pi +\frac { \pi} {2\beta}$, as long as $\beta$ is a positive  irrational number.

If $\beta>0$ is a rational number, then we write $\beta=p/q$ with $p$ and $q$ being relatively prime. The change of variable $t^{1/q}=s$ gives
\begin{equation}\label{G-rational}
  G_\beta(z)=\int^\infty_0 \frac {q s^{q-1} e^{-s^p}}{s^q+\zeta^q} d\zeta, ~~\zeta^q=z.
\end{equation}Let $\omega_l=e^{(2l+1)\pi i/q}$, $l=0,1,\cdots, q-1$. Since
\begin{equation*}
  s^q+\zeta^q=(s-\omega_0\zeta)   (s-\omega_1\zeta)\cdots (s-\omega_{q-1}\zeta),
\end{equation*}taking logarithmic derivative yields
\begin{equation*}
  \frac {q s^{q-1}}{s^q+\zeta^q}=\sum^{q-1}_{l=0} \frac 1 {s-\omega_l\zeta}.
\end{equation*}Thus, $G_\beta(z)$ in \eqref{G-rational} can be written as
\begin{equation*}
  G_\beta(z)=\sum^{q-1}_{l=0} G_p\left ( e^{(2l+1-q)\pi i/q} \zeta\right ).
\end{equation*}The problem is now reduced to consider the function $G_p(z)$, where $p$ is a positive integer. In this case, the integral for $G_p(z)$ can again be broken  into two parts, one over the interval $(0, 1)$ and the other on the infinite interval $(1, \infty)$. The integral on  $(1, \infty)$ can be asymptotically evaluated as before. For the integral on $(0,1)$, we expand the exponential function $e^{-t^p}$ into a power series. This leads to infinitely  many integrals of the form
\begin{equation*}
  I_n(z)= \int^1_0\frac {t^{n p}}{t+z} dt,~~n=0,1,2,\cdots.
\end{equation*} For each $n\geq 1$, we have
\begin{equation*}
  \frac 1 {t+z}=\frac 1 t-\frac z {t^2}+\cdots+ (-1)^{np-1}\frac {z^{np-1}}{t^{np}}+(-1)^{np}\frac {z^{np}}{t^{np}(t+z)},
\end{equation*}and hence
\begin{equation*}
  I_n(z)=\sum^{np-1}_{k=0} \frac {(-1)^kz^k}{np-k}  +(-1)^{np}z^{np}\left [ \ln (1+z)-\ln z\right ].
\end{equation*}Asymptotic expansion of $G_p(z)$ can thus be obtained for small values of $z$.

\begin{figure}[h]
\begin{center}
\includegraphics[width=7.7cm]{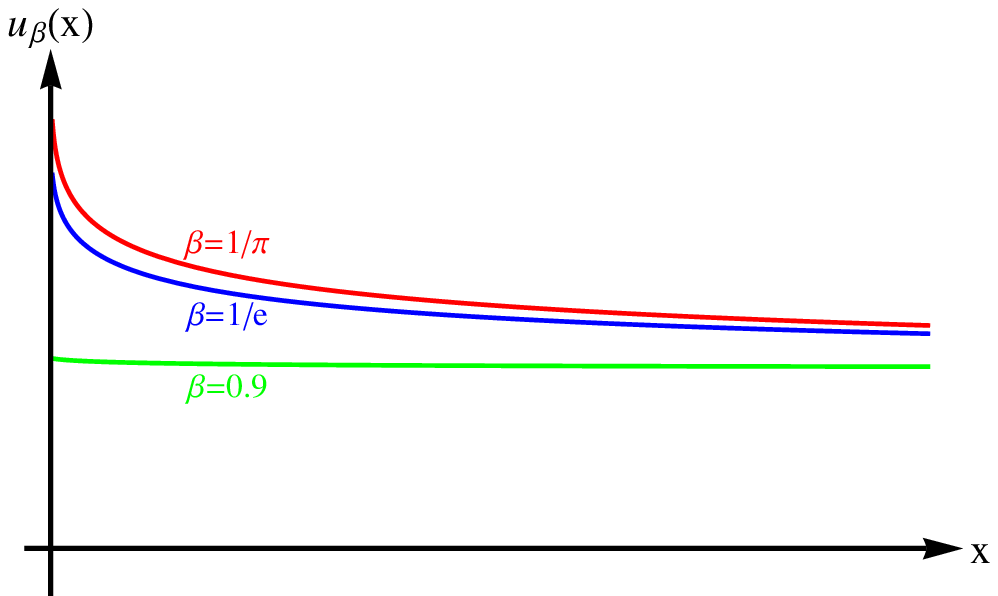}\hskip .5cm \includegraphics[width=7.7cm]{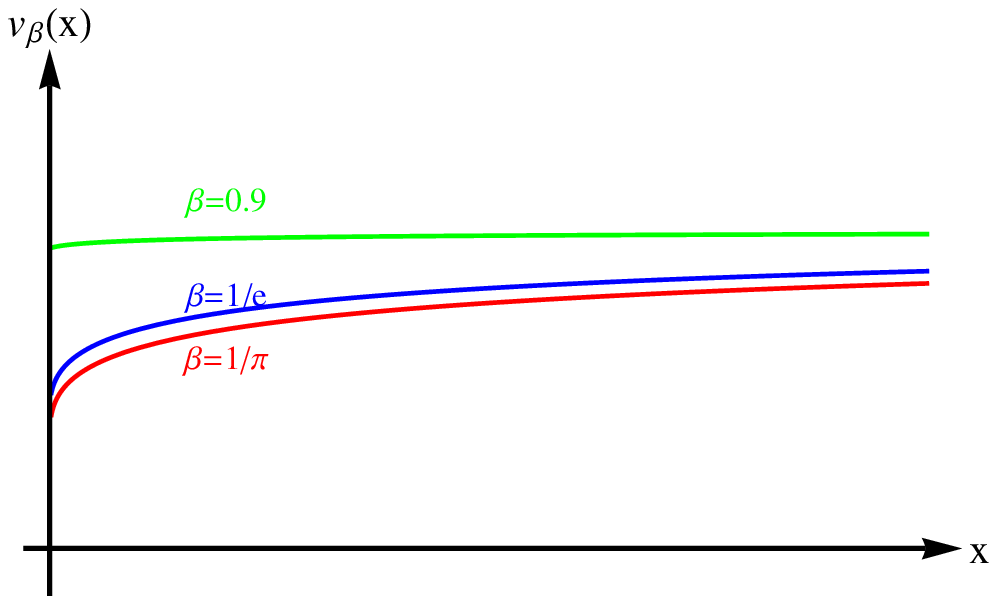}
 \caption{The solution profiles of $u_\beta(x)$ (left) and $v_\beta(x)$ (right).}
 \label{figures-for-u-v}
\end{center}
\end{figure}

\noindent
\begin{que}
We have used numerical simulation to draw the graphs of $u_\beta(x)$ for various values $\beta\in (0, 1)$; see Figure \ref{figures-for-u-v}. From the graphs, it is reasonable to conjecture that the function $u_\beta(x)$ is monotonically decreasing in $x\in [0, \infty)$, and also that the initial value $u_\beta(0)$ is decreasing in $\beta\in (0, 1)$. What we need now are proofs of these statements. From Figure \ref{figures-for-u-v}, corresponding conjectures can be made for the function $v_\beta(x)$.
\end{que}

\section*{Acknowledgements}
We would like to thank Dr. Xiang-Sheng Wang for drawing the graphs and helping us with the numerical simulation.   The research of Yu-Qiu Zhao was
  supported in part by the National
Natural Science Foundation of China under grant numbers
10871212 and  11571375.

\end{document}